\documentclass[12pt]{amsart}
\usepackage{amsmath,amsthm,amsfonts,amssymb}
\usepackage{graphicx,psfrag}
\begin{document} 
\newcommand{\B}{{\mathbb B}}
\newcommand{\C}{{\mathbb C}}
\newcommand{\N}{{\mathbb N}}
\newcommand{\Q}{{\mathbb Q}}
\newcommand{\Z}{{\mathbb Z}}
\renewcommand{\P}{{\mathbb P}}
\newcommand{\R}{{\mathbb R}}
\newcommand{\tensor}{\otimes}
\newcommand{\rc}{\subset}
\newcommand{\rank}{\mathop{rank}}
\newcommand{\trace}{\mathop{tr}}
\newcommand{\dimc}{\mathop{dim}_{\C}}
\newcommand{\Lie}{\mathop{Lie}}
\newcommand{\Auto}{\mathop{{\rm Aut}_{\mathcal O}}}
\newcommand{\Aut}{\mathop{{\rm Aut}}}
\newcommand{\alg}[1]{{\mathbf #1}}
\newtheorem*{definition}{Definition}
\newtheorem{claim}{Claim}
\newtheorem{corollary}{Corollary}
\newtheorem*{conjecture}{Conjecture}
\newtheorem*{SpecAss}{Special Assumptions}
\newtheorem{example}{Example}
\newtheorem*{remark}{Remark}
\newtheorem*{observation}{Observation}
\newtheorem*{fact}{Fact}
\newtheorem*{remarks}{Remarks}
\newtheorem{lemma}{Lemma}
\newtheorem{proposition}{Proposition}
\newtheorem{theorem}{Theorem}
\title{%
An Example related to Brody's theorem
}
\author {J\"org Winkelmann}
\begin{abstract}
We discuss an example related to the method of Brody.
\end{abstract}
%
\address{%
J\"org Winkelmann \\
 Institut Elie Cartan (Math\'ematiques)\\
 Universit\'e Henri Poincar\'e Nancy 1\\
 B.P. 239\\
 F-54506 Vand\oe uvre-les-Nancy Cedex\\
 France
}
\email{jwinkel@member.ams.org\newline\indent{\itshape Webpage: }%
http://www.math.unibas.ch/\~{ }winkel/
}
\thanks{
{\em Acknowledgement.}
The author wants to thank 
V.~Bangert and B.~Siebert for
the invitation to the workshop in Freiburg
in September 2003.
}
\maketitle
\section{Introduction}
\subsection{Bloch principle}
In one-dimensional function theory there is a general philosophy
which supposedly goes back to A.~Bloch (see e.g.~\cite{Z},\cite{Bl}): 
If $P$ is a sufficiently
reasonable class of
holomorphic maps or functions, then the following statements
should be equivalent:
\begin{enumerate}
\item
Every map in class $P$ defined on the complex line $\C$ is constant.
\item
The set of all maps in class $P$ defined on the unit disk
$\Delta=\{z\in\C:|z|<1\}$ is a {\em normal family}.
\end{enumerate}
(A family of maps is called a ``normal family'' if every sequence
in it is either compactly divergent or contains a subsequence
which converges uniformly on compact sets. A sequence of maps
$f_n:X\to Y$ between topological spaces is ``compactly divergent'', 
if for every pair of compact subsets $K\subset X$, $C\subset Y$ there
are only finitely many $f_n$ with $f_n(K)\cap C\ne\{\}$.)

For example, every bounded holomorphic function on $\C$ is constant
by Liouville's theorem and due to Montel's theorem the family of all
bounded holomorphic functions  on $\Delta$ is a normal family.
Thus the Bloch principle is valid for the family $P$ of all bounded
holomorphic functions with values in $\C$.

\subsection{Brody's theorem}
Let $Y$ be a complex manifold. 
It is called ``taut'' if the family
of all holomorphic maps $f:\Delta\to Y$ is a
normal family. Let us from now on assume that $Y$ is compact.
Then being ``taut'' is easily seen to be equivalent
with hyperbolicity in the sense of Kobayashi.
The theorem of Brody (see \cite{B}) states that this is furthermore
equivalent with
the property that every holomorphic map from $\C$ to $Y$ is constant.
In other words: Brody's theorem states that the Bloch principle hold
for the class of holomorphic maps with values in a (fixed) compact complex
manifold $Y$.

Now we may raise the question: What about holomorphic maps to a compact
complex manifold fixing some given base points?
Given a compact complex manifold $Y$ and a point $y\in Y$, 
let us consider the following two statements:
\begin{itemize}
\item
Every holomorphic map $f:\C\to Y$ with $f(0)=y$ is constant.
\item
The family of all holomorphic maps $f:\Delta\to Y$
with $f(0)=y$ is a normal family.
\end{itemize}

Are they equivalent?

Using the notion of the infinitesimal Koba\-yashi-Royden pseudometric
as introduced in \cite{R} this can be reformulated into the following
question:
{\em
``If the infinitesimal Koba\-yashi-Royden peusdometric on a
compact complex manifold $Y$
degenerates for some point $y\in Y$, does this imply that there
exists a holomorphic map $f:\C\to Y$ with $y\in f(\C)$?''}

Thanks to Brody's
theorem it is clear
 that there exists some non-constant holomorphic map $f:\C\to Y$
if the Koba\-yashi-Royden pseudometric is degenerate at some point $y$
of $Y$.
But it is not clear that $f$ can be chosen in such a way that $y$ is in the
image or at least in the closure of the image.
Of course, at first it looks absurd that degeneracy of the Koba\-yashi-Royden
pseudometric at one point $y$ should only imply the existence of
a non-constant holomorphic map to some part of $Y$ far away of $y$ and
should not imply the existence of a non-constant map $f:\C\to Y$ whose
image comes close to $y$.

Thus one is led to postulate
\begin{conjecture}
Let $X$ be a compact complex manifold, $x\in X$.
Assume that the infinitesimal Koba\-yashi-Royden pseudometric
is degenerate on $T_xX$.

Then there exists a non-constant holomorphic map $f:\C\to X$
with $f(0)=x$.
\end{conjecture}

\subsection{Bounded derivatives}
Let $X$ be a complex manifold equipped with a hermitian metric $h$.
For each holomorphic map $f:\C\to X$ and each point $z\in\C$ we may
now calculate the norm of the derivatie $Df$ at $z$ with respect
to the euclidean metric on $\C$ and $h$ on $X$.
Let $P$ be a class of holomorphic maps $f:\C\to (X,h)$ with
bounded derivatives (i.e.~for every $f\in P$ there is a number $C>0$ such
that the inequality $||Df_z||<C$ holds for all $z\in\C$).
Let $f:\C\to X$ be a non-constant map
in this class $P$. Via $f_n(x)=f(nx)$ this map $f$ yields a non-normal
family of maps $f_n:\Delta\to X$.

Now let $P'$ denote the set of those maps in $P$ for which the derivative
(calculated with respect to the euclidean metric on $\C$ resp.~$\Delta$
and the hermitian metric on $X$) is bounded.
For each of the $f_n$ defined above the derivative is clearly bounded,
since $\Delta$ is relatively compact in $\C$, and $f_n:\C\to X$
extends through the boundary. Thus $f_n$ is a non-normal family in $P'$.
If the Bloch principle holds for $P'$, this implies the existence
of a non-constant holomorphic map $F:\C\to X$ in $P$.

Thus: {\em If the Bloch principle holds for $P'$, the existence of
a non-constant holomorphic map $f$ in $P$ implies the existence
of a non-constant holomorphic map $F$ in $P$ with the additional
property that $||DF||$ is bounded.}

Brody's theorem implies that this is indeed true if, given a compact
complex manifold $X$, we consider the set $P$ of all holomorphic
maps with values in $X$.

However, we will give an example of a compact complex manifold
$X$, an open subset $\Omega$ and a point $x\in\Omega$ such that this
property does not hold if $P$ is chosen as the family of all holomorphic
maps $f$ with image contained in $\Omega$ and $f(0)=x$.

\subsection{Reparametrization}
The key method for proving a Bloch principle is the following:
Let $f_n:\Delta\to Y$ be a non-normal family. Then
we look for an increasing sequence of disk $\Delta_{r_n}$
which exhausts $\C$ (i.e. $\lim{r_n}=+\infty$) and
a sequence of holomorphic maps $\alpha_n:\Delta_{r_n}\to\Delta$
such that a subsequence of $f_n\circ\alpha_n$ converges
(locally uniformly) to a non-constant holomorphic map from
$\C$ to $Y$.

For the proof of his theorem Brody used this idea, taking combinations
of affine-linear maps with automorphisms of the disk for the $\alpha_n$.

Zalcman (\cite{Z})
investigated other reparametrizations where the $\alpha_n$
themselves are affin-linear maps, a concept which has the advantage
that it can also be applied to harmonic maps.

\subsection{Subvarieties of abelian varieties}
Let $A$ be a complex abelian variety (i.e. a compact complex torus
which is simultaneously a projective algebraic variety)
and $X$ a subvariety. Let $E$ denote the union of all translates
of complex subtori of $A$ which are contained in $X$.
It is known that this union is either all of $X$ or a proper
algebraic subvariety (\cite{Kw}).

Since $A$ is a compact complex torus there is a flat hermitian metric
on $A$ induced by the euclidean metric on $\C^g$ via $A\simeq\C^g/\Gamma$.
A holomorphic map $f:\C\to A$ has bounded derivative with respect to
this metric if and only if it is induced by an affine-linear map
from $\C$ to $\C^g$.

From this, one can deduce that $f(\C)\subset E$ for every holomorphic
map $f:\C\to X$ with bounded derivative. Given the previous
considerations about the Bloch principle, it is thus natural
to conjecture:
\begin{conjecture}
For every non-constant holomorphic map $f:\C\to X$ the image is contained
in $E$. The Koba\-yashi-pseudodistance on $X$ is a distance outside $E$.
\end{conjecture}

For example, this statement is a consequence of the more general
conjecture VIII.I.4 by S. Lang in \cite{LNT}.
In the context of classification theory the above statement has also
be conjectured by F. Campana (\cite{C},\S 9.3).

In the spirit of the analogue between diophantine geometry and
entire holomorphic curves as pointed out by Vojta
\cite{V}, the
conjecture above is also supported by the famous result of Faltings
(\cite{F})
with which he solved the Mordell conjecture. This result states the
following: If we assume that $A$ und $X$ are defined over a number field
$K$, then with only finitely many exceptions every $K$-rational point
of $X$ is contained in $E$.

\subsection{Our example}
We construct an example of the following type:
There is a compact complex manifold $X$, equipped with some hermitian
metric, an open subset $\Omega$ and a point $p\in\Omega$.
There exists a non-constant holomorphic map $f:\C\to\Omega$ with
$f(0)=p$. Via $f_n(z)=f(nz)$ this yields a non-normal family
of holomorphic maps $f_n:\Delta\to\Omega$ with bounded derivatives
such that $f_n(0)=p$.

But there is no non-constant holomorphic map $f:\C\to\Omega$ with
$f(0)=p$ {\em and bounded derivative}.

\section{The example}
\subsection{Statement of main results}
We construct an example which shows that Brody reparametrization
sometimes necessarily changes the image of the curve.

\begin{theorem}
There exists a compact complex hermitian manifold $(T,h)$ 
and open subsets
$\Omega_2\subset\Omega_1\subset T$ such that:
\begin{enumerate}
\item
$\Omega_2$ is not dense in $\Omega_1$
and neither is $\Omega_1$ in $T$.
\item
For every point $p\in\Omega_1$
there is a non-constant holomorphic map
$f:\C\to\Omega_1$ with $p=f(0)$.
\item
If $f:\C\to T$ is a non-constant
holomorphic map with bounded derivative
(with respect to the euclidean metric on $\C$ and $h$ on $T$) and
$f(\C)\subset\bar\Omega_1$,
then $f(\C)\subset\bar\Omega_2$.
\end{enumerate}
\end{theorem}

Recall that Brody's method, starting from any holomorphic
map from $\C$ to $T$,
yields a holomorphic map from $\C$ to $T$ with bounded
derivative. 
Thus this examples provides a picture in which Brody's method
really changes the properties of $f:\C\to T$ fundamentally.

Responding to some additional questions which may be asked, we prove
a little bit more.

\begin{theorem}\label{thm-detailed}
There exists a compact complex torus $T$, equipped with
a flat hermitian metric $h$ 
and open subsets
$\Omega_2\subset\Omega_1\subset T$ such that:
\begin{enumerate}
\item
$\Omega_2$ is not dense in $\Omega_1$
and neither is $\Omega_1$ in $T$.
\item
For every point $p\in\Omega_1$ and every $v\in T_p\Omega_1$
there is a non-constant holomorphic map
$f:\C\to\Omega_1$ with $p=f(0)$, $v=f'(0)$ and $\bar\Omega_1=\overline{f(\C)}$.
\item
If $f:\C\to T$ is a non-constant
holomorphic map with bounded derivative
(with respect to the euclidean metric on $\C$ and $h$ on $T$) and
$f(\C)\subset\bar\Omega_1$,
then $f(\C)\subset\bar\Omega_2$.
Moreover $f$ is affine-linear and $\overline{f(\C)}$ is a 
closed analytic subset of $T$.
\end{enumerate}
\end{theorem}

We remark that this implies in particular that the infinitesimal
Koba\-yashi-Royden pseudometric vanishes identically on $\Omega_1$.

Furthermore, it provides examples of holomorphic maps from $\C$
into a compact complex torus with a rather ``bad'' image:
The closure of the image with respect to the euclidean topology
is $\overline{\Omega_2}$ and thus a set with non-empty interior
such that the complement has also non-empty interior.
This is in strong contrast to the Zariski-analytic closure:
By the theorem of Green-Bloch-Ochiai for every holomorphic
map $f$ from $\C$ to a compact complex torus $T$ the closure
of the image $f(\C)$
with respect to the analytic Zariski topology (i.e.~the smallest
closed analytic subset of $T$ containing $f(\C)$)
is always a translated subtorus of $T$.

We will now describe our example.

We precede the construction with some elementary observations 
about tori:
Let $T=\C^n/\Lambda$ be a torus, equipped with the flat euclidean
metric and the corresponding distance function $d_T(\ ,\ )$.
Let 
\[
\rho=\frac{1}{2}\min_{\gamma\in\Lambda\setminus\{0\}}||\gamma||.
\]
This is the {\sl injectivity radius}, in other words $\rho$ is
the largest real number such that
the natural projection $\pi:\C^n\to T$ induces a homeomorphism
between the ball
\[
B_\epsilon(\C^n;0)=\{v\in\C^n:||v||<\epsilon\}
\]
and
\[
B_\epsilon(T;e)=\{x\in T:d_T(x,e)<\epsilon\}
\]
for all $\epsilon<\rho$.
Evidently, the injectivity radius $\rho$ is a lower bound
for the {\em diameter}
\[
\rho\le diam=\max_{x,y\in T}d_T(x,y)
\]
If we pass from $T$ to a subtorus $S\subset T$, the injectivity radius
can only increase, while the diameter can only decrease.
As a consequence we obtain:
\begin{lemma}
Let $T$ be a compact (real or complex) torus
with injectivity radius
$\rho$. 
Then for every real 
positive-dimensional
subtorus $S\subset T$ the diameter
\[
diam(S)=\max_{x,y\in S}d_T(x,y)
\]
is at least $\rho$.

Furthermore, if
$0<\epsilon<\rho$ and $x\in T$, then the ball $B_\epsilon(T;x)$
contains no translate of any positive-dimensional real subtorus of $T$.
\end{lemma}

Before giving the details of the construction
of our example, let us try to express its idea in a drawing:

\psfrag{Ea}{$E''$}
\psfrag{Eb}{$E'$}
\psfrag{W}{$W$}
\psfrag{pa}{${\leftarrow}$}
\psfrag{pb}{$\downarrow$}
\psfrag{Oa}{$\Omega_1$}
\psfrag{Ob}{$\Omega_2$}
\psfrag{Oc}{$\Sigma$}

\smallskip
\includegraphics{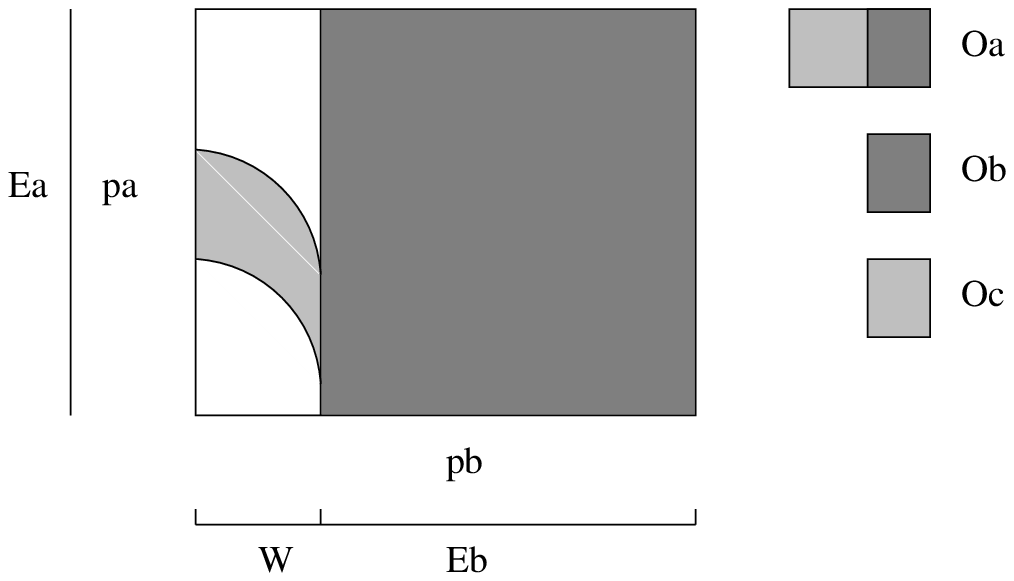}

\smallskip\par
Now let us start the precise construction of the example.
Let $E'=\C/\Gamma'$ and $E''=\C/\Gamma''$ 
be elliptic curves and $T=E'\times E''$.
Let $\pi':\C\to E'$, $\pi'':\C\to E''$ 
and $\pi=(\pi',\pi''):\C^2\to T$
denote the natural projections.
We assume that $E'$ is not isogenous to $E''$. 
(For example, we might choose $E'=\C/\Z[i]$ and
$E''=\C/\Z[\sqrt 2i]$.)
Then $E'\times\{0\}$ and $\{0\}\times E''$ are the only non-trivial
complex subtori of $T$.

Now $T=\C^2/\Gamma$ with $\Gamma=\Gamma'\times\Gamma''$.
The compact complex torus $T$ carries a  hermitian metric $h$
induced by the euclidean metric on $\C^2$ 
(i.e.~$h=dz_1\tensor d\bar z_1 +dz_2\tensor d\bar z_2$).
The associated distance function is called $d$, the
injectivity radius $\rho$ is defined as explained above.

We choose numbers $0<\rho'<\rho''<\rho$ and define
$W=B_{\rho'}(E',e)$.

Furthermore we choose $0<\delta<\frac{1}{2}\rho$ and we choose
a holomorphic map $\sigma:\C\to E''$ such that
there exist complex numbers $t,t'\in B_{\rho'}(\C,0)$ 
(i.e.~$|t|,|t'|<\rho'$)
and 
\[
d_{E''}(\sigma(t),\sigma(t'))>2\delta.
\]
(This is possible since $2\delta$ is smaller than the injectivity
radius $\rho$ of $T$ which in turn is a lower bound for the
diameter of $E''$).

We denote by $s:\C\to \C$ a holomorphic function such that
$\sigma=\pi''\circ s$.

Since $\pi':\C\to E'$ restricts to an isomorphism between
$B_{\rho}(\C,0)$ and $B_{\rho}(E',e)$, the holomorphic maps
$s$ and $\sigma$ induce maps from $B_\rho(E',e)$ to $\C$ resp.~$E''$.
By abuse of notation these maps will also be denoted by $s$
resp.~$\sigma$.

Now define $\Omega_2=(E'\setminus\bar W)\times E''$ and 
$\Omega_1=\Omega_2\cup\Sigma$
with
\[
\Sigma=\{(x,y):x\in \bar W, y\in E'', d_{E''}(y-\sigma(x))<\delta\}
\]

Let us now fix some point $p\in\Omega_1$
and $v=(v_1,v_2)\in T_p(T)=\C^2$. 
We have to show that there exists a holomorphic map $f$
as stipulated in $(2)$ of theorem~\ref{thm-detailed}.

Let $(p_1,p_2)\in\C^2$ be a point mapped on $p$
by $\pi:\C^2\to T$. If $p\in\Sigma$, we require $|p_1|\le\rho'$
and $|s(p_1)-p_2|<\delta$ and define $\delta'=\delta-|s(p_1)-p_2|$.
If $p\not\in\Sigma$, we require $|p_1|>\rho''$ and define $\delta'=\delta$.

As the next step, we will choose a pair of entire functions
$(Q,H)$.

\begin{claim}
There is a pair of entire functions $(Q,H)$ with the following properties:
\begin{enumerate}
\item
$Q$ is a non-constant polynomial,
\item
$(Q(0),H(0))=(p_1,p_2)$ and
\item
$(Q'(0),H'(0))$ and $v$ are parallel.
\item
If $p\in\Sigma$, we require furthermore 
that $(Q(z),H(z)+y)\in\pi^{-1}(\Sigma)$
for all $z$ and $y$ with $|Q(z)|\le\rho'$ and $|y|\le\frac{1}{2}\delta'$.
\end{enumerate}
\end{claim}

Let us first discuss the case where $p\not\in\Sigma$.

Then it suffices to choose
\[
Q(z)= z^2+v_1z+p_1
\]
and
\[
H(z)=v_2z+p_2.
\]

If $p\in\Sigma$, we proceed as follows:
First, for $r,t\in \C$ we define 
\[
Q_t(z)=(z+t)^2+p_1-t^2
\]
and
\[
H_{r,t}(z)= p_2-s(p_1)+s(Q_t(z))+rz.
\]

We will set $Q=Q_t$ and $H=H_{r,t}$ for appropriately chosen
parameters $r,t$.

Evidently $Q_t$ is a polynomial for any choice of $t$.
Furthermore $(Q_t(0),H_{r,t}(0))=(p_1,p_2)$ independent of the choice
of $r,t$:

\[
Q_t(0)=t^2+p_1-t^2=p_1
\]
and
\[
H_{r,t}(0)= p_2-s(p_1)+s(p_1)+0 =p_2.
\]
Let $\Phi_{r,t}=(Q_t,H_{r,t})$.
We have
\[
\Phi'_{r,t}(0)=(Q_t'(0),s'(Q_t(0))Q_t'(0)+r)=(2t,2s'(p_1)t+r)
\]
Observe that
\[
(r,t)\mapsto \frac{2t}{2s'(p_1)t+r}
\]
defines a meromorphic function on $\C^2$ with a point of
indeterminacy at $(0,0)$.
This is true regardless of the value of $s'(p_1)$.

Thus every neighborhood of $(0,0)$ contains
a point $(r,t)\ne(0,0)$ such that $\Phi'_{r,t}(0)$ is a non-zero
multiple of $v$. 

Next we note that $(t,z)\mapsto Q_t(z)$ defines a proper map
from $\overline{B_1(\C,0)}
\times\C$ to $\C$. Therefore there is a constant $C>0$
such that $|z|<C$, whenever there exists a parameter $t$ such that
$|t|\le 1$ and $|Q_t(z)|\le\rho$.

It is therefore possible to choose two numbers $r,t$ in such a way
that
\begin{enumerate}
\item $\Phi'_{r,t}(0)$ is a non-zero multiple of $v$,
\item $|t|<1$ and
\item $|2rC|<\delta'$.
\end{enumerate}

Now assume that $z,y\in\C$ with $|Q_t(z)|\le\rho'$
and $|y|<\frac{1}{2}\delta'$.
By the definition of the constant $C$, this implies $|z|<C$.
Let $(w_1,w_2)=\Phi_{r,t}(z)+(0,y)$.
Then
\begin{multline*}
|w_2-s(w_1)|=|p_2-s(p_1)+rz+y| 
< |p_2-s(p_1)| + |rC| + \frac{1}{2}\delta'
<\\
< (\delta-\delta') + \frac{1}{2}\delta' + \frac{1}{2}\delta'=\delta
.
\end{multline*}

Now $|w_2-s(w_1)|<\delta$ in combination with $|w_1|=|Q_t(z)|\le\rho'$
implies $\pi(w_1,w_2)\in\Sigma$.
Hence $\Phi_{r,t}(z)+(0,y)\in\pi^{-1}(\Sigma)$ 
under this assumption.
Thus the claim is proved.
Q.E.D.

Our next step is to construct a closed subset $A$ of $\C$ to which
we will apply Arakelyan approximation.

Let $A_0$ be the union of $\overline{B_{\rho''}(0)}$
and $\overline{B_{\rho'}(\gamma)}$ for all $\gamma\in\Gamma'$.
If $p\not\in\Sigma$, then $p_1\not\in A_0$.
Hence in this case we can choose $\eta>0$ such that 
$\overline{B_\eta(p_1)}$ is disjoint to $A_0$ and define $A_1$
as the union of $A_0$ with this closed ball $\overline{B_\eta(p_1)}$.
If $p\in\Sigma$, we simply take $A_1=A_0$.

Next we choose dense countable subsets $S_1\subset int(\Sigma)$ 
(where $int(\Sigma)$ denotes the interior of $\Sigma$) and
$S_2\subset\Omega_2$. We observe that $\C\setminus A_1$ projects
surjectively onto $E'\setminus\overline{W}$ and that the fibers
of this projection are infinite discrete subsets of $\C$.
For this reason we can find sequences $a_n,b_n$ in $\C$
such that
\[
S_2=\{\pi(a_n,b_n):n\in\N\}
\]
and all the $a_n$ are distinct elements of $\C\setminus A_1$
with $\lim_{n\to\infty}|a_n|=+\infty$.
It follows that
\[
\Theta=\{a_n:n\in\N\}
\]
is a discrete subset of $\C$ which has empty intersection with $A_1$.
We define $A_2=A_1\cup\Theta$.

We fix a bijection $\xi:\Gamma'\setminus\{0\}\stackrel{\sim}\rightarrow
S_1$ and an enumeration $n\mapsto \gamma_n$ of $\Gamma'\setminus\{0\}$.
Then we can choose sequences of complex numbers $c_n,d_n$ such that
the following properties hold  for all $n\in\N$
\begin{enumerate}
\item $\pi(c_n,d_n)=\xi(\gamma_n)$,
\item $|c_n-\gamma_n|<\rho'$ and
\item $|d_n-s(c_n)|<\delta$.
\end{enumerate}

We define $A=Q^{-1}(A_2)$.
\begin{claim}
Arakelyan approximation is applicable to $A$,
i.e. $\{\infty\}\cup(\C\setminus A)$ is connected
and locally connected.
\end{claim}

Observing that we can deform $B_{\rho''}(\C,0)$ to $B_{\rho'}(\C,0)$,
we deduce from prop.~\ref{arak-pol} that $Q^{-1}(A_0)$ has the
desired property. Now $A$ and $Q^{-1}(A_0)$ differ only by removing
the preimage of a closed disc and by removing a discrete countable
set (namely $Q^{-1}(\Theta)$). This can not destroy connectivity,
hence not only $\{\infty\}\cup(\C\setminus Q^{-1}(A_0))$
but also $\{\infty\}\cup(\C\setminus A)$ is connected and locally
connected.
Thus the claim is proved.

We will now define a continuous function $h$ on $A$, which is
holomorphic in its interior, and which we will then approximate
by an entire function, using Arakelyan's theorem.

If $p\not\in\Sigma$,
we take $h(z)=H(z)$
on $Q^{-1}(\overline{B_\eta(p_1)})$
and $h=s$ on  $Q^{-1}(\overline{B_{\rho''}(0)})$.

If $p\in\Sigma$, we define $h$ on $Q^{-1}(\overline{B_{\rho''}(0)})$
as $H(z)$.

Next, for every $n\in\N$, we define $h(z)$ as
\[
h(z)=s(Q(z)-\gamma_n)+d_n-s(c_n)
\]
whenever $|Q(z)-\gamma_n|\le\rho'$.

Finally, we define $h$ on $Q^{-1}(\Theta)$ by stipulating that
$h(z)=b_n$ whenever $Q(z)=a_n$ for a number $n\in\N$.

By the construction of $(Q,H)$
we know that  $\pi(Q(0),h(0))=p$ and that $(Q'(0),h'(0))$ is a multiple of
$v$. The choice of $h$ implies moreover that $S_1\cup S_2$ is contained
in the image of $z\mapsto \pi(Q(z),h(z))$.

Next we define 
a continuous positive function $\epsilon:A\to\R^+$ as follows:

\begin{itemize}
\item
$\epsilon\equiv 1$ on $Q^{-1}(\overline{B_\eta(p_1)})$ if $p\not\in\Sigma$.
\item
$\epsilon\equiv\frac{1}{2}\delta'$ on $Q^{-1}(\overline{B_{\rho''}(0)})$.
\item
$\epsilon(z)=\frac{1}{n}$ if $Q(z)=a_n$.
\item
$
\epsilon(z)=\min\left\{\frac{1}{n},\frac{1}{2}\left(\delta-|d_n-s(c_n)|\right)
\right\}
$
whenever $|Q(z)-\gamma_n|\le\rho'$.
\end{itemize}

Using prop.~\ref{arak-inter}, we deduce that there exists an entire
function $F:\C\to\C$ such that
\begin{enumerate}
\item
$|F(z)-h(z)|<\epsilon(z)$ for all $z\in A$.
\item
$F(0)=h(0)$ and $F'(0)=h'(0)$.
\end{enumerate}

By the second condition we obtain that $\pi(Q(0),F(0))=p$ and that
$(Q'(0),F'(0))$ is a multiple of $v$.
The first condition ensures that $\pi(Q(z),F(z))\in\Omega$
for all $z\in\C$. It also ensure that the image is dense:
Indeed, let $w\in\Omega_2$. Then there is a sequence of points
in $S_2$ converging to $w$. But $S_2=\{\pi(a_n,b_n):n\in\N\}$
and the construction of $F$ implies that for every $n\in\N$ there
exists a number $z_n\in\C$ such that $Q(z_n)=a_n$ and 
$|F(z_n)-b_n|<\frac{1}{n}$. It follows that there is a subsequence
$z_{n_k}$ such that $\lim_k \pi(Q(z_{n_k}),F(z_{n_k}))=w$.
If $w\in\Sigma$, we argue similarily, with $S_1$ in the 
role of $S_2$.
Thus the whole set $\Omega_1$ is in the closure of the image
of the map $z\mapsto \pi(Q(z),F(z))$ from $\C$ to $T$.

Finally, let $\mu$ be a complex number such that
$\mu(Q'(0),F'(0))= v$ and define
\[
f(z)= \pi\left(Q(\mu z),F(\mu z)\right)
\]
Then $f:\C\to\Omega_1$ is a holomorphic map with the desired
properties.

\section{Arakelyan Approximation with interpolation}
We will need a slight improvement of Arakelyan's theorem.

We recall the theorem of Arakelyan (see \cite{A}):
\begin{theorem}
Let $A$ be a closed subset of $\C$, $U=\P_1(\C)\setminus A$,
$\epsilon:A\to\R^+$ a continuous function and $f_0:A\to \C$
a continuous function which is holomorphic in the interior
of $A$.
Assume that $U$ is connected and locally connected.

Then there exists a holomorphic function $F:\C\to\C$ with
 $|F(z)-f(z)|<\epsilon(z)$ for all $z\in A$.
\end{theorem}

We want to verify that Arakelyan's theorem is applicable in our
situation.
\begin{proposition}\label{arak-pol}
Let $\Gamma$ be a lattice in $\C$ and $\rho'$ a real number with
\[
0<\rho' < \rho=\frac{1}{2}\min_{\gamma\in\Gamma\setminus\{0\}}
|\gamma|
\]
Let $A'=\{z\in\C: d(z,\Gamma)\le \rho'\}=\cup_{\gamma\in\Gamma}
\overline{B_{\rho'}(\C,\gamma)}$, $P:\C\to\C$ a non-constant
polynomial and $U=\{\infty\}\cup (\C\setminus P^{-1}(A'))$.

Then $U$ is connected and locally connected.
\end{proposition}
\begin{proof}
First we want to verify that $U$ contains no bounded connected component.
Indeed, assume that there is such a connected component $C$.
Its boundary $\partial C$ is a connected set mapped into
\[
\cup_{\gamma\in\Gamma}\overline{B_{\rho'}(\gamma)}
\]
by P. 
This is a disjoint union due to the choice of $\rho'$.
Hence continuity of $P$ implies that there is one element $\gamma\in\Gamma$
such that $|P(z)-\gamma|\le \rho'$ for all $z\in \partial C$.

But $C\subset U$ implies $|P(z)-\gamma|>\rho'$ for all $z\in C$, $\gamma\in
\Gamma$.
This is in contradiction with the maximum principle
for the holomorphic function $P$.
Hence there can not exist a bounded connected component $C\subset U$.

\psfrag{A1}{$A'$}
\psfrag{B1}{$B_n(\C,0)$}
\psfrag{V1}{$R_n$}
\includegraphics{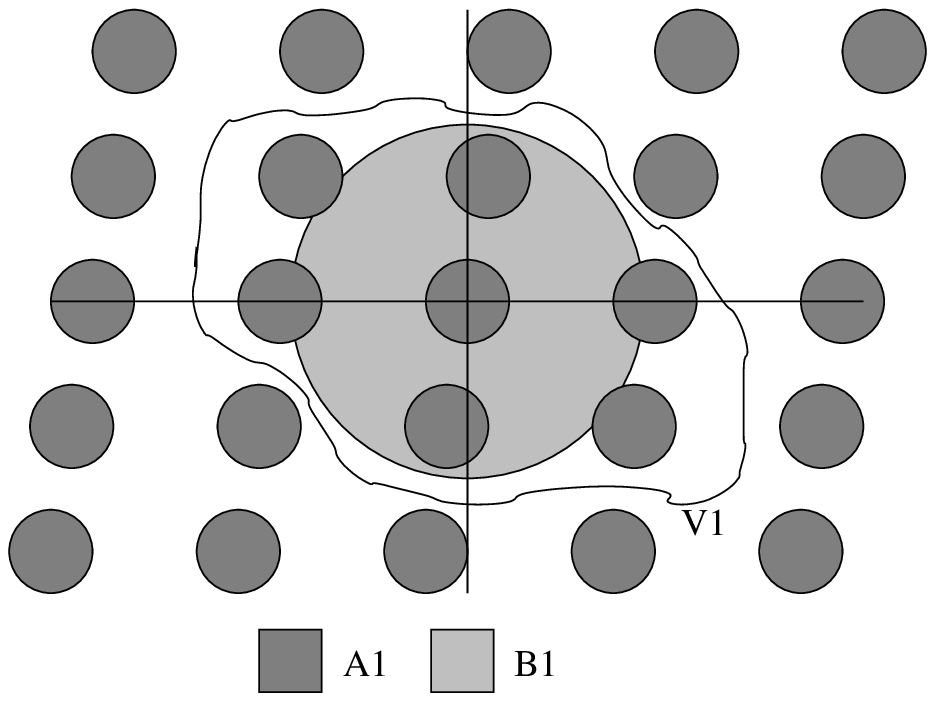}

For each $n\in\N$ we 
choose a simple closed curve $R_n\subset\C\setminus A'$
such that the open bounded subset $V_n\subset\C$
which is enclosed by $R_n$ has the property that $B_n(\C,0)\subset V_n$.

The ramification locus 
\[
Z=\{z\in\C:\exists w\in\C : P(w)=z, P'(w)=0\}
\]
is
a finite set. Let $N_0=\max\{|z|:z\in Z\}$.
Now the restriction of $P$ to $P^{-1}(\C\setminus V_n)\to \C\setminus V_n$
is an unramified covering of degree $d=\deg(P)$ for all $n>N_0$.
As a polynomial map, $P$ extends to a proper map $\bar P:\C\cup\{\infty\}
\to\C\cup\{\infty\}$. For a suitably chosen local coordinate $w$ at
$\infty$ the map $P$ near $\infty$ can be described as $w\mapsto w^d$.
Using this fact and the fact that by construction each curve $R_n$
defines a generator for 
\[
\pi_1(\C^*)\simeq\pi_1(\C\setminus\overline{B_n(\C,0)})
\]
we can conclude that $P^{-1}(R_n)$ is connected for all $n>N_0$.
Then $P^{-1}(\Omega)$ is connected for every open subset $\Omega\subset
\C\setminus\overline{B_n(\C,0)}$ with $R_n\subset\Omega$.

In particular 
\[
W_n=U\setminus P^{-1}(V_n)
\]
is connected for all $n>N_0$.
The collection of all these open sets $W_n$ constitutes an neighborhood
basis of $U$ at $\infty$, implying that $U$ is locally connected
at infinity.
Furthermore, the connectedness of the sets $W_n$ implies that $U$ is only
one unbounded connected component. Since we have already seen that $U$
is no bounded connected component, this completes the proof that $U$
is connected and locally connected.
\end{proof}

\begin{proposition}\label{arak-inter}
Let $A$ be a closed subset in $\C$, $A\ne\C$, 
and suppose that for every function
$f$ on $A$ which is holomorphic in its interior and every continuous
map $\epsilon:A\to\R^+$ there is an entire function $F:\C\to\C$
with $|F(z)-f(z)|<\epsilon(z)$ for all $z\in A$.

Let $q$ be a point in the interior of $A$. Then we can find such an
entire function $F$ with the additional properties $F(q)=f(q)$
and $F'(q)=f'(q)$.
\end{proposition}

\begin{proof}
Let $U=\{\infty\}\cup\left(\C\setminus A\right)$.
By assumption $U$ is connected and locally connected at infinity.
Let $p\in\C\setminus A$ and 
let $W$ be a bounded connected 
open subset of $\C$ containing both $p$ and $q$.
Choose $\delta>0$ such that 
\[
\delta<\min\left\{ d(q,\partial W),d(q,\partial A),d(p,q)\right\}
\]
and define
\[
\tilde A=\overline{B_{\delta}(q)}\cup A\setminus W
\]
and
\[
\tilde U=\{\infty\}\cup\left(\C\setminus\tilde A\right)
=U\cup\left(W\setminus \overline{B_{\delta}(q)}\right).
\]
Now both $U$ and $\left(W\setminus \overline{B_{\delta}(q)}\right)$
are connected, and their intersection is non-empty, since it contains $p$.
Therefore
$\tilde U$ is connected.
Moreover, $\tilde U$ is locally connected
at infinity, because it coincides with $U$ near $\infty$.
Thus we have Arakelyan approximation for $\tilde A$.

\psfrag{qa}{$q$}
\psfrag{pa}{$p$}
\psfrag{A}{$A$}
\psfrag{W}{$W$}
\includegraphics{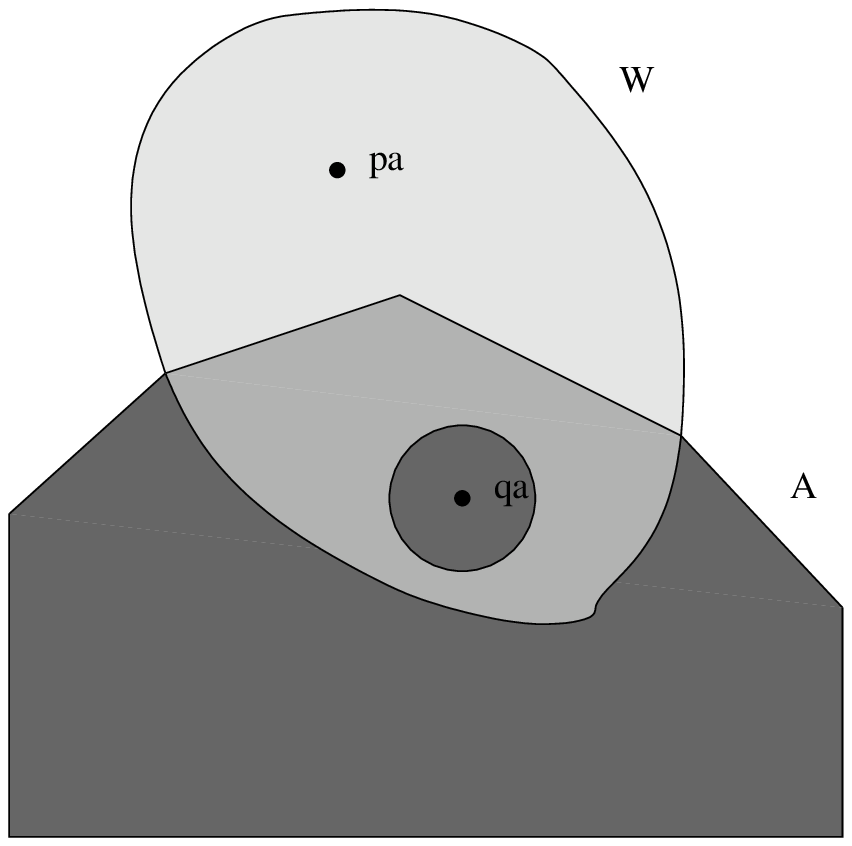}

We choose constants $\xi_0,\xi_1\in\C\setminus\{0\}$ such that
\[
|\xi_0|<\frac{1}{16}\epsilon(z)
\]
and
\[
|\xi_1(z-q)|<\frac{1}{16}\epsilon(z)
\]
for all $z\in \overline{B_{\delta}(q)}$.

Then we define functions $g,h:\tilde A\to\C$
via
\[
g(z)=
\begin{cases} \xi_0 & \text{ if }z\in \overline{B_{\delta}(q)}\\
0 & \text{ else } 
\end{cases}
\]
and
\[
h(z)=
\begin{cases} \xi_1(z-q) & \text{ if }z\in \overline{B_{\delta}(q)}\\
0 & \text{ else } 
\end{cases}
\]

Clearly, $g$ and $h$ are continuous and holomorphic in the interior of
$\tilde A$.
The choice of $\xi_0,\xi_1$ implies that $|g(z)|<\frac{1}{16}\epsilon(z)$
and $|h(z)|<\frac{1}{16}\epsilon(z)$ for all $z\in A$.

By the Arakelyan property we find sequences of entire functions $g_n,h_n:\C\to\C$ such that
\[
|g_n(z)-g(z)|<\frac{1}{8n}\epsilon(z)
\]
and
\[
|h_n(z)-h(z)|<\frac{1}{8n}\epsilon(z)
\]
for all $n\in\N$, $z\in\tilde A$.
Locally uniform convergence on $\tilde A$ implies that inside the interior of 
$\tilde A$ the derivatives converge as well.
Hence we obtain
\[
\lim_{n\to\infty}\begin{pmatrix} g_n(q) & h_n(q) \\ g_n'(q) & h_n'(q)
		 \end{pmatrix}
=\begin{pmatrix} g(q) & h(q) \\ g'(q) & h'(q) 
 \end{pmatrix}
=\begin{pmatrix} \xi_0 & 0 \\ 0 & \xi_1
 \end{pmatrix}.
\]
Thus, for $n$ sufficiently large the vectors $(g_n(q),g_n'(q))$
and $(h_n(q),h_n'(q))$ are linearly independent.

Next we observe that $A\setminus\tilde A$ is relatively compact 
in $A$. Therefore, for  
sufficiently large numbers $n,C$ 
the functions $\alpha=\frac{1}{C}g_n$
and $\beta=\frac{1}{C}h_n$ have the following properties:

\begin{enumerate}
\item
$\alpha,\beta$ are entire functions,
\item
$|\alpha(z)|,|\beta(z)|<\frac{1}{8}\epsilon(z)$ for all $z\in A$,
and
\item
the vectors $(\alpha(q),\alpha'(q))$ and $(\beta(q),\beta'(q))$
are linearly independent.
\end{enumerate}

By the approximation property for $A$ there are 
sequences of entire functions $\alpha_n,\beta_n,f_n:\C\to\C$
such that
\[
\max\{ |\alpha_n(z)-\alpha(z)|, |\beta_n(z)-\beta(z)|,
|f_n(z)-f(z)|\}<\frac{1}{n}\epsilon(z)
\]
for all $n\in\N$, $z\in A$.
The locally uniform convergence of $\lim\alpha_n=\alpha$,
$\lim\beta_n=\beta$ and $\lim f_n=f$ on $A$ implies that in the
interior of $A$ the respective derivatives converge as well.
In particular, this happens at $q$.
Hence the matrix
\[
A_n=\begin{pmatrix}
\alpha_n(q) & \beta_n(q) \\ \alpha_n'(q) & \beta_n'(q)
\\
       \end{pmatrix}
\]
converges to
\[
\lim_{n\to\infty}A_n=
A=\begin{pmatrix}
\alpha(q) & \beta(q) \\ \alpha'(q) & \beta'(q)
\\
       \end{pmatrix}.
\]
Since $A$ is invertible, it follows that
$A_n$ is likewise invertible for all sufficiently large $n$.
Hence we can define 
(for sufficiently large $n$)
sequences $\lambda_n,\mu_n$ via
\[
\begin{pmatrix} \lambda_n \\ \mu_n \end{pmatrix}
=
A_n^{-1}\cdot
\begin{pmatrix} f(q)-f_n(q) \\ f'(q)-f'_n(q) \end{pmatrix}.
\]
Now $\lim f_n=f$, $\lim f_n'=f'$ and $\lim A_n^{-1}=A^{-1}$.
Therefore $\lim\lambda_n=0=\lim\mu_n$.

Thus we can choose a natural number $N\in\N$ with the following
properties:

\begin{enumerate}
\item $A_N$ is invertible,
\item $|\lambda_N|,|\mu_N|<1$,
\item and $N>4$.
\end{enumerate}

We define 
\[
F(z)=f_N(z)+\lambda_N\alpha_N(z)+\mu_N\beta_N(z).
\]

By the choice of $\lambda_n,\mu_n$ we have 
\begin{multline*}
\begin{pmatrix}
F(q) \\ F'(q) 
\end{pmatrix}
=\begin{pmatrix}
f_N(q)+\lambda_N\alpha_N(q)+\mu_N\beta(q) \\
f'_N(q)+\lambda_N\alpha'_N(q)+\mu'_N\beta'(q)
 \end{pmatrix}
=\\
=\begin{pmatrix}
f_N(q) \\ f_N'(q) 
\end{pmatrix}
+
A_N\cdot \begin{pmatrix} \lambda_N \\ \mu_N \end{pmatrix}
= 
\begin{pmatrix}
f(q) \\ f'(q) 
\end{pmatrix}.
\end{multline*}

Furthermore 
\begin{multline*}
|F(z)-f(z)| \le |f_N(z)-f(z)| + 
|\lambda_N|\cdot\left(|\alpha_N(z)-\alpha(z)|+|\alpha(z)|\right)+\\
+|\mu_N|\cdot\left(|\beta_N(z)-\beta(z)|+|\beta(z)|\right)\\
\le \frac{1}{N}\epsilon(z) +  \frac{1}{N}\epsilon(z) +  \frac{1}{8}\epsilon(z) 
+ \frac{1}{N}\epsilon(z)+ \frac{1}{8}\epsilon(z)
<\left(\frac{3}{4}+\frac{2}{8}\right)\epsilon(z)=\epsilon(z)
\end{multline*}
for all $z\in A$.
Thus $F$ is an entire function with the desired properties.
\end{proof}

\end{document}